\title{Reciprocal cyclotomic polynomials}
\author{Pieter Moree}
\documentclass[12pt]{article}
\usepackage{amssymb, latexsym, amsfonts}
\def\@ptsize{2}
\newtheorem{Thm}{Theorem}

\newtheorem{Lem}{Lemma}

\newtheorem{cor}{Corollary}

\newcommand{\qed}{\hfill $\Box$}

\begin{document}
\date{}
\maketitle
{\def\thefootnote{}
\footnote{{\it Mathematics Subject Classification (2000)}.
11B83, 11C08}}

\begin{abstract}
\noindent Let $\Psi_n(x)$ be the monic polynomial having precisely all
non-primitive $n$th roots of unity as its simple zeros. One
has $\Psi_n(x)=(x^n-1)/\Phi_n(x)$, with $\Phi_n(x)$ the $n$th cyclotomic polynomial.
The coefficients of $\Psi_n(x)$ are integers that like the coefficients of $\Phi_n(x)$ tend to
be surprisingly small in absolute value, e.g. for $n<561$ all coefficients of $\Psi_n(x)$
are $\le 1$ in absolute value. We establish various properties of the
coefficients of $\Psi_n(x)$.
\end{abstract}
\section{Introduction}
The $n$th cyclotomic polynomial $\Phi_n(x)$ is defined by
$$\Phi_n(x)=\prod_{1\le j\le n\atop (j,n)=1}(x-\zeta_n^j)=
\sum_{k=0}^{\varphi(n)}a_n(k)x^k,$$
where $\varphi$ is Euler's totient function
and $\zeta_n$ a primitive $n$th root of unity. 
The coefficients $a_n(k)$
are known to be integral. The study of the $a_n(k)$ began with the
startling observation that for small $n$ we have $|a_n(k)|\le 1$. The first
counterexample to this inequality occurs for $n=105$: $a_{105}(7)=-2$.
The amazement over the smallness of $a_n(m)$ was eloquently phrased by D. Lehmer \cite{DL}
who wrote: `The smallness of $a_n(m)$ would appear to be one of the fundamental conspiracies
of the primitive $n$th roots of unity. When one considers that $a_n(m)$ is a sum of
$({\phi(n)\atop m})$ unit vectors (for example $73629072$ in the case of $n=105$ and
$m=7$) one realizes the extent of the cancellation that takes place'.\\
\indent We define
$\Psi_n(x)$ by
$$\Psi_n(x)=\prod_{1\le j\le n\atop (j,n)>1}(x-\zeta_n^j)=
\sum_{k=0}^{n-\varphi(n)}c_n(k)x^k.$$
Note that $\Psi_n(x)=(x^n-1)/\Phi_n(x)$. The identity $x^n-1=\prod_{d|n}\Phi_d(x)$ shows that
$$\Psi_n(x)=\prod_{d|n, d<n}\Phi_d(x),$$
and thus the coefficients of $\Psi_n(x)$ are integers.\\ 
\indent Note that for $|x|<1$ we have
$${1\over \Phi_n(x)}=-\Psi_n(x)(1+x^n+x^{2n}+\cdots ).$$
Since $n>n-\varphi(n)$, it follows that the Taylor coefficients of $1/\Phi_n(x)$ are
periodic with a period dividing $n$. This allows one to easily reformulate the results
on the coefficients of $\Psi_n(x)$ obtained in this paper to the Taylor coefficients
of $1/\Phi_n(x)$ as well.\\
\indent The purpose of
this note is to show that the non-primitive roots, like the primitive ones, conspire and study
the extent to which this is the case.
\section{Some basics}
Note that
\begin{equation}
\label{productje}
x^n-1=\prod_{d|n}\prod_{1\le j\le n\atop (j,n)=d}(x-\zeta_n^j)=\prod_{d|n}\Phi_{n\over d}(x)=
\prod_{d|n}\Phi_d(x).
\end{equation}
It follows from this identity that
\begin{equation}
\label{opzijnkop}
\Psi_n(x)={x^n-1\over \Phi_n(x)}=\prod_{d|n\atop d<n}\Phi_d(x).
\end{equation}
We infer that $\Psi_n(x)\in \mathbb Z[x]$. 
\begin{Lem}
\label{beginnie}
Let $n>1$. We have
$$\Psi_n(x)=-\prod_{d|n\atop d<n}(1-x^d)^{-\mu({n\over d})}.$$
\end{Lem}
{\it Proof}.
By applying M\"obius inversion one deduces from (\ref{productje}) that
\begin{equation}
\label{basiccyclo}
\Phi_n(x)=\prod_{d|n}(x^d-1)^{\mu({n\over d})}.
\end{equation}
On using that $\sum_{d|n}\mu(n/d)=0$, we infer that
$\Phi_n(x)=\prod_{d|n}(1-x^d)^{\mu({n\over d})}$, from
which the result follows on invoking (\ref{opzijnkop}). \qed\\

\noindent Let rad$(n)=\prod_{p|n}p$ be the {\it radical} of $n$. {}From the
previous lemma it is not difficult to arrive at
the next result, see e.g. Thangadurai \cite{Thanga} for 
the proof of the corresponding result for $\Phi_n(x)$.
\begin{Lem}
\label{blup}
Let $n>1$. We have:\\
{\rm 1)} $\Psi_{2n}(x)=(1-x^n)\Psi_n(-x)$ if $n$ is odd;\\
{\rm 2)} $\Psi_{pn}(x)=\Psi_n(x^p)$ if $p|n$;\\
{\rm 3)} $\Psi_{pn}(x)=\Psi_n(x^p)\Phi_n(x)$ if $p\nmid n$;\\
{\rm 4)} $\Psi_n(x)=\Psi_{{\rm rad}(n)}(x^{n\over {\rm rad}(n)})$;\\
{\rm 5)} $\Psi_n(x)=-\Psi_n({1\over x})x^{n-\varphi(n)}$.
\end{Lem} 
Put $V_n=\{c_n(k):0\le k\le n-\varphi(n)\}$. If $n>1$ then by part 5 of the latter
lemma we have that $a\in V_n$ implies that $-a\in V_n$. It also gives that if
$n-\varphi(n)$ is even, then $c_n((n-\varphi(n))/2)=0$.
\begin{Lem}
\label{opstart}
If $n=1$, then $V_n=\{1\}$. If $n$ is a prime, then  $V_n=\{-1,1\}$.
In the remaining cases we have $\{-1,0,1\}\subseteq V_n$.
\end{Lem}
{\it Proof}. If $n=1$, then $\Psi_n(x)=1$. If $n$ is a prime, then
$\Psi_n(x)=x-1$. Next assume that $n$ has at least
two (not necessarily distinct) prime divisors. Note that this implies
that $n-\varphi(n)\ge 2$. Note that $\Psi_n(x)$ is monic and that
 $\Psi_n(0)=-1$ by Lemma \ref{beginnie}. It thus remains to be shown
that $0\in V_n$. In case $n$ is not squarefree we have
$\Psi_n(x)=-1+O(x^2)$ by Lemma \ref{beginnie} and thus $c_n(1)=0$. If $n$
is odd and squarefree and $\mu(n)=-1$, then by Lemma \ref{beginnie} we
find $\Psi_n(x)=-1+x+O(x^3)$ and hence $c_n(2)=0$ (here we use
that $n-\varphi(n)\ge 2$). If $n$
is odd and squarefree and $\mu(n)=1$, then by Lemma \ref{beginnie} we
find $$\Psi_n(x)={(x^p-1)\over 1-x}(1+O(x^{p+1})),$$
where $p$ is the smallest prime divisor of $n$ and hence $c_n(p)=0$. 
Since $p\le n-\varphi(n)$ it follows that $0\in V_n$.
In
case $n$ is even and squarefree we invoke part 1 of Lemma \ref{blup}
to complete the proof. \qed\\

\noindent It is not difficult to prove that, as $x$ tends to infinty,
$$\sum_{n\le x}{\varphi(n)\over n}\sim {x\over \zeta(2)}=x{6\over \pi^2}.$$
Thus the average degree of $\Phi_n(x)$ and $\Psi_n(x)$ is ${6\over \pi^2}n$, 
respectively $(1-{6\over \pi^2})n$.
We have ${6\over \pi^2}=0.60792710\cdots$ and
$1-{6\over \pi^2}=0.3920728\cdots$.

\subsection{ (Reciprocal) cyclotomic polynomials of low order}
We define the order of $\Phi_n(x)$ and
$\Psi_n(x)$ to be the number, $\omega_1(n)$, of distinct odd prime divisors of $n$.
Instead of saying that $f$ has order 3, we sometimes say that $f$ is ternary.
We define the {\it height} of a polynomial $f$ in $\mathbb Z[x]$, $h(f)$,  to be the
maximum  absolute value of the coefficients of $f$.
In case $h(f)=1$ we say that $f$ is {\it flat}.\\ 
\indent Low order examples (in the remainder of this section $p<q<r$ 
will be primes):\\
$\Psi_1(x)=1$;\\
$\Psi_p(x)=-1+x$;\\
$\Psi_{pq}(x)=-1-x-x^2-\ldots-x^{p-1}+x^q+x^{q+1}+\ldots+x^{p+q-1}.$\\
These examples in combination with parts 1 and 4 of Lemma \ref{blup} establish
the correctness of the following result.
\begin{Lem}
\label{kleiner3}
If $\Psi_n(x)$ is of order $\le 2$, then $\Psi_n(x)$ is flat.
\end{Lem}
We like to point out that $\Psi_{pq}(x)$ has a rather simpler structure than
$\Phi_{pq}(x)$. 
It can be shown, see e.g. Carlitz \cite{C}, Lam and Leung \cite{LL} and Thangadurai \cite{Thanga}, that
\begin{eqnarray}
\Phi_{pq}(x)&=&\sum^{\rho}_{i=0}x^{ip}\sum^{\sigma}_{j=0}x^{jq}-x^{-pq}
\sum^{q-1}_{i=\rho+1}x^{ip}\sum^{p-1}_{j=\sigma+1}x^{jq}\nonumber\cr
&=&\sum^{\rho}_{i=0}x^{ip}\sum^{\sigma}_{j=0}x^{jq}-x\sum^{q-2-\rho}_{i=0}x^{ip}\sum^{p-2-\sigma}_{j=0}x^{jq},
\end{eqnarray}
where $\rho$ and $\sigma$ are the unique nonnegative integers for 
which $(p-1)(q-1)=\rho p+\sigma q$ (note that $\rho\le q-2$ and $\sigma \le p-2$). 
As a consequence we have
the following evaluation of the coefficients $a_{pq}(k)$.
\begin{Lem}
\label{phiflat}
Let $p<q$ be odd primes. Let $\rho$ and $\sigma$ be the unique nonnegative integers for 
which $(p-1)(q-1)=\rho p+\sigma q$.
Then 
$$a_{pq}(k)=\cases{1 & if $k=ip+jq$ for some $0\le i\le\rho,~0\le j\le \sigma$;\cr
-1 & if $k=1+ip+jq$ for some $0\le i\le q-2-\rho,0\le j\le p-2-\sigma$;\cr
0 & otherwise.}$$
\end{Lem}
Using the latter lemma it is easy to show that if
$\Phi_n(x)$ is of order $\le 2$, then $\Phi_n(x)$ is flat.\\
\indent {}For the convenience of the reader we will prove that there are unique non-negative
integers such that $(p-1)(q-1)=\rho p+\sigma q$ (this
proof is taken from Ram\'irez Alfons\'in's book
\cite[p. 34]{RA}, with the observation that in case $p$ and $q$ are primes
the auxiliary polynomial $Q(x)$ equals $\Phi_{pq}(x)$). We let $r(n)$ be the number of 
representations of $n$ in the form $n=px+qy$ with $x,y\ge 0$. We have
$$R(x)=\sum_{i=0}^{\infty}r(i)x^i={1\over (1-x^p)(1-x^q)}.$$
Note that $R(x)(x^{pq}-1)(x-1)=\Phi_{pq}(x)$. By L'H\^opital's rule we find that
$\Phi_{pq}(1)=1$ and hence we have that
$${\Phi_{pq}(x)-1\over x-1}=\sum_{i=0}^{pq-p-q}g(i)x^i,$$
with $g(pq-p-q)=1$. On the other hand,
\begin{eqnarray}
{\Phi_{pq}(x)-1\over x-1}&=&R(x)(x^{pq}-1)+{1\over 1-x};\nonumber\cr
&=&\sum_{i=0}^{\infty}r(i)x^{pq+i}+\sum_{i=0}^{\infty}(1-r(i))x^i;\nonumber\cr
&=&\sum_{i=0}^{pq-1}(1-r(i))x^i+\sum_{i=pq}^{\infty}(r(i-pq)+1-r(i))x^i.\nonumber
\end{eqnarray}
On comparing the two expressions for $(\Phi_{pq}(x)-1)/(x-1)$ we arrive at various
conclusions. {}First of all we see that $r((p-1)(q-1))=1$. Secondly it allows one to
compute the {\it Frobenius number} $g(p,q)$. Given relatively prime positive integers
$a_1,\ldots,a_n$ the largest natural number that is not representable as a non-negative
integer combination of $a_1,\ldots,a_n$ is called the Frobenius number and denoted
by $g(a_1,\ldots,a_n)$. On noting that $r(i-pq)\le r(i)$ comparison of the two expressions
for $(\Phi_{pq}(x)-1)/(x-1)$ shows that $r(pq-p-q)=0$ and $r(pq-p-q+i)\ge 1$ for
$i\ge 1$, which yields $g(p,q)=pq-p-q$.\\
\indent By 
Lemma \ref{beginnie} we have $$\Psi_{pqr}(x)={(x-1)(1-x^{pq})(1-x^{pr})(1-x^{qr})
\over (1-x^p)(1-x^q)(1-x^r)}.$$
This can be written as
\begin{equation}
\label{pqr}
\Psi_{pqr}(x)=(x-1)\Big(\sum_{j_1=0}^{q-1}x^{j_1p}\Big)\Big(\sum_{j_2=0}^{r-1}x^{j_2q}\Big)
\Big(\sum_{j_3=0}^{p-1}x^{j_3r}\Big).
\end{equation}
Alternatively we can write, by part 3 of Lemma \ref{blup},
\begin{equation}
\label{pqr2}
\Psi_{pqr}(x)=\Phi_{pq}(x)\Psi_{pq}(x^r).
\end{equation}
Let the {\it denumerant} be defined as the number of non-negative integer representations
of $m$ by $a_1,a_2,\ldots,a_n$. Denote it by $d(m;a_1,\ldots,a_n)$. For $m<pq$ we
infer from (\ref{pqr}) that $c_{pqr}(k)=d(m-1;p,q)-d(m;p,q)$. For more on denumerants
see Chapter 4 of Ram\'irez Alfons\'in \cite{RA}.
\begin{Lem}
\label{flauw}
Let $p<q<r$ be odd primes. If $0\le k<r$, then we have
$c_{pqr}(k)=-a_{pq}(k)\in \{-1,0,1\}$.
\end{Lem}
{\it Proof}. Immediate from (\ref{pqr2}), $\Psi_{pq}(0)=-1$  and Lemma \ref{phiflat}. \qed\\

\noindent The following result also relates $c_{pqr}(k)$ to $a_{pq}(k)$ in case $k>r$.
(If $k$ is outside the range $[0,\ldots,\varphi(n)]$ respectively, $[0,\ldots,n-\varphi(n)]$,
then we put $a_n(k)=0$, respectively $c_n(k)=0$.
\begin{Lem}
\label{verbinding}
Let $p<q<r$ be odd primes. Put $\tau=(p-1)(r+q-1)$. Suppose that $qr>\tau$.
If $k\le \tau$, then
$$c_{pqr}(k)=-\sum_{j=0}^m a_{pq}(k-jr),$$
with $m$ the unique integer such that $mr\le k<(m+1)r$. Furthermore,
\begin{equation}
\label{prop}
c_{pqr}(\tau-k)=c_{pqr}(k)
\end{equation}
and $c_{pqr}(k+qr)=-c_{pqr}(k)$.
If $\tau<k<qr$, then $c_{pqr}(k)=0$.
\end{Lem}
{\it Proof}. We have
\begin{equation}
\label{glik2}
\Psi_{pqr}(x)=\Phi_{pq}(x)(1+x^r+\ldots+x^{(p-1)r})(x^{qr}-1).
\end{equation} 
Write 
\begin{equation}
\label{boe2}
\Phi_{pq}(x)(1+x^r+\dots+x^{(p-1)r})=\sum_{k=0}^{\tau}e_{pqr}(k)x^k.
\end{equation}
Note that the polyniomal in (\ref{boe2}) of degree $\tau$ and selfreciprocal.
If $k\le \tau$, then $c_{pqr}(k)=-e_{pqr}(k)$ and $c_{pqr}(k+qr)=e_{pqr}(k)$.
On combining all these observations the result easily follows. \qed\\

\noindent In 1895 Bang \cite{Bang} proved that $h(\Phi_{pqr}(x))\le p-1$. The
same bound applies to the height of $\Psi_{pqr}(x)$.
\begin{Thm}
\label{upper}
The height of $\Psi_{pqr}(x)$ is at most $p-1$. More precisely, we have
$$h(\Psi_{pqr}(x))\le\Big[{(p-1)(q-1)\over r}\Big]+1.$$
\end{Thm}
{\it Proof}. By (\ref{pqr2}) we find that
\begin{equation}
\label{startie}
c_{pqr}(k)=\sum_{j=0}^{[k/r]}a_{pq}(k-jr)c_{pq}(j).
\end{equation}
The number of $j$ for which $0\le k-jr\le \varphi(pq)$ is
$$\le \Big[{\varphi(pq)\over r}\Big]+1=\Big[{(p-1)(q-1)\over r}\Big]+1\le p-2+1=p-1.$$
The proof is finished since $|a_{pq}(k-jr)|\le 1$ by Lemma \ref{phiflat} and 
$|c_{pq}(j)|\le 1$ by the identity
$\Psi_{pq}(x)=-1-x-x^2-\ldots-x^{p-1}+x^q+x^{q+1}+\ldots+x^{p+q-1}$.\qed\\

\noindent We have seen that on average the degree of $\Phi_n(x)$ is less than that 
of $\Psi_n(x)$. It is left to the reader to show that if $p<q<r$ are odd primes,
then deg$(\Psi_{pqr}(x))<{\rm deg}(\Phi_{pqr}(x))$, except when
$pqr\in  \{105,165,195\}$.

\section{Beiter's conjecture and its reciprocal analogue}
In 1971 Sister Marion Beiter \cite{Beiter} put forward the conjecture that
if $p<q<r$ are odd primes, then $\Phi_{pqr}(x)$ is of
height at most $(p+1)/2$. As she pointed out, her conjecture is
true for $p\le 5$. She also showed that the height is $\le p-\lfloor p/4 \rfloor$.
Bachman \cite{B1} showed that if either $q$ or $r$ is congruent to $\pm 1$ or $\pm 2$ modulo $p$,
then the height is $\le (p+1)/2$. 
H. M\"oller \cite{HM2} gave explicit examples of polynomials $\Phi_{pqr}(x)$, for
every $p$, with a prescribed coefficient equal to $(p+1)/2$. This shows that
the conjecture is best possible, if true. 
More precisely, M\"oller showed that if $q\equiv -2({\rm mod~}p)$, $r\equiv
-(p-1)(q-1)/2 ({\rm mod~}pq)$, then $a_{pqr}((p-1)(qr+1)/2)=(p+1)/2$.
For further results and references see Bachman 
\cite{B1, B2}. In general Beiter's conjecture remains unresolved.\\
\indent The following result gives the analogue of the Beiter conjecture for
the reciprocal polynomials.
\begin{Thm}
\label{main}
Let $p<q<r$ be odd primes. Then
$h(\Psi_{pqr}(x))=p-1$ iff 
$$q\equiv r\equiv \pm 1({\rm mod~}p){\rm ~and~}r<{(p-1)\over (p-2)}(q-1).$$
In the remaining cases $h(\Psi_{pqr}(x))<p-1$.
\end{Thm} 
\begin{cor} Suppose that $h(\Psi_{pqr}(x))=p-1$ and $q+2p$ is a prime, then
also $h(\Psi_{pq(q+2p)}(x))=p-1$.
\end{cor}
By the above theorem and Dirichlet's theorem on arithmetic progressions it follows
that for every prime $p\ge 3$ there are infinitely many pairs $(q,r)$ such that
$h(\Psi_{pqr}(x))=p-1$.\\
\indent Theorem \ref{main} follows from two theorems that deal with the
necessity, respectively sufficiency part of its iff statement in combination with Theorem
\ref{upper}.
\begin{Thm}
\label{ext}
If $h(\Psi_{pqr}(x))=p-1$, then 
$$q\equiv r\equiv \pm 1({\rm mod~}p){\rm ~and~}r<{(p-1)\over (p-2)}(q-1).$$
\end{Thm}
{\it Proof}. Let $j_{\min}$ be the smallest $j$ such that $k-jr\le \varphi(pq)$ and
$j_{\max}$ be the largest $j$ such that $k-jr\ge 0$. Then we can write (\ref{startie})
as $$c_{pqr}(k)=\sum_{j=j_{\min}}^{j_{\max}}a_{pq}(k-jr)c_{pq}(j).$$
From $k-j_{\max}r\ge 0$ and $k-j_{\min}r\le (p-1)(q-1)$ we infer
that $(j_{\max}-j_{\min})r\le (p-1)(q-1)<(p-1)r$ and hence 
$j_{\max}-j_{\min}\le p-2$. In order to have $c_{pqr}(k)=p-1$ for some $k$ we must 
have $j_{\max}-j_{\min}=p-2$. Thus $(j_{\max}-j_{\min})r=(p-2)r\le (p-1)(q-1)$.
Since $(p-2)r$ is odd and $(p-1)(q-1)$ is even it follows that
$$r<{(p-1)\over (p-2)}(q-1).$$
Let $k$ be such that $|c_{pqr}(k)|=p-1$. Then we must have that
$c_{pq}(j)\ne 0$ for $j_{\min}\le j\le j_{\max}$.
It follows from this that 
the pair $(j_{\min},j_{\max})$ must be one of the following: 
$(0,p-2),~(1,p-1),~(q,q+p-2),~(q+1,q+p-1)$, and that 
$c_{pq}(j_{\min})=c_{pq}(j_{\min}+1)=\ldots=c_{pq}(j_{\max})$. Thus we have
$$p-1=|c_{pqr}(k)|=\Big|\sum_{j=j_{\min}}^{j_{\max}}a_{pq}(k-jr)\Big|.$$
We now make a case distinction according to whether $a_{pq}(k-jr)=1$ for 
$j_{\min}\le j\le j_{\max}$, or $a_{pq}(k-jr)=-1$ for every
$j_{\min}\le j\le j_{\max}$.\\
{\tt {}First case}. {}For every $j_{\min}\le j\le j_{\max}$ we have $a_{pq}(k-jr)=1$.\\
By Lemma \ref{phiflat} it follows that there must be non-negative integers $i_m$ and $j_m$ with $0\le i_m\le \rho$
and $0\le j_m\le \sigma$ such that
$$\cases{k-j_{\max}r&=$i_1p+j_1q$;\cr  k-(j_{\max}-1)r&=$i_2p+j_2q$; \cr \cdots &= $\cdots$ \cr
k-j_{\min}r&=$i_{p-1}p+j_{p-1}q,$}$$
Now if we would have $j_{m_1}=j_{m_2}$ for $m_1\ne m_2$ by subtracting the corresponding
equations we infer that $p|r$, a contradiction. Thus we must have 
$\{j_1,\ldots,j_{p-1}\}=\{0,1,\ldots,p-2\}$ and hence $\sigma=p-2$. It follows that
$q\equiv -1({\rm mod~}p)$ and $\rho=(q-p+1)/p$. Now select $m_1$ and $m_2$ such
that $j_{m_2}=j_{m_1}+1$. On substracting the corresponding equations we infer that
$\alpha r=\beta p+q$ for some integers $\alpha$ and $\beta$ with $-\rho \le \beta\le \rho$.
Note that $p-1\le \beta p+q<2q-p+1<2r$. It follows that $\alpha=1$
and $r=\beta p+q$ and hence $r\equiv q\equiv -1({\rm mod~}p)$.\\
{\tt Second case}. {}For every $j_{\min}\le j\le j_{\max}$ we have $a_{pq}(k-jr)=-1$.\\
By Lemma \ref{phiflat} it then
follows that there must be non-negative integers $i_m$ and $j_m$ with $0\le i_m\le q-2-\rho$
and $0\le j_m\le p-2-\sigma$ such that
$$\cases{k-j_{\max}r&=$1+i_1p+j_1q$;\cr  k-(j_{\max}-1)r&=$1+i_2p+j_2q$; \cr \cdots &= $\cdots$ \cr
k-j_{\min}r&=$1+i_{p-1}p+j_{p-1}q,$}$$
For the same reason as above we must have 
$\{j_1,\ldots,j_{p-1}\}=\{0,1,\ldots,p-2\}$. This implies $\sigma=0$. It follows
that $q\equiv 1({\rm mod~}p)$ and $\rho=(p-1)(q-1)/p$ and thus $\rho':=q-2-\rho=(q-p-1)/p$.
Now select $m_1$ and $m_2$ such
that $j_{m_2}=j_{m_1}+1$. On substracting the corresponding equations we infer that
$\alpha r=\beta p+q$ for some integers $\alpha$ and $\beta$ with $-\rho' \le \beta\le \rho'$.
Note that $p+1\le \beta p+q<2q-p-1<2r$. It follows that $\alpha=1$
and $r=\beta p+q$ and hence $r\equiv q\equiv 1({\rm mod~}p)$.\qed

\begin{Thm}
\label{extreme}
Let $p<q<r$ be odd primes such 
that $r<(p-1)(q-1)/(p-2)$. If
$q\equiv -1({\rm mod~}p)$ and $r\equiv -1({\rm mod~}p)$, then
$$c_{pqr}(k)=\cases{-1-m & for $0\le m\le p-2$, $k=mr$;\cr
0 & for $k=2$;\cr 
m+1 & for $0\le m\le p-2$, $k=(m+q)r$,}$$
and $V_{pqr}=\{-(p-1),-(p-2),\ldots,p-2,p-1\}$.\\
\indent If $q\equiv 1({\rm mod~}p)$ and $r\equiv 1({\rm mod~}p)$, then
$$c_{pqr}(k)=\cases{1+m & for $0\le m\le p-2$, $k=1+mr$;\cr
0 & for $k=2$;\cr 
-1-m & for $0\le m\le p-2$, $k=1+(m+q)r$,}$$
and $V_{pqr}=\{-(p-1),-(p-2),\ldots,p-2,p-1\}$.
\indent 
\end{Thm} 
{\it Proof}. {}From the proof of Lemma \ref{opstart} it follows that $c_{pqr}(2)=0$.\\
{\tt {}First case}. Assume that $q\equiv r\equiv -1({\rm mod~}p)$.\\
Note that $\rho=(q-p+1)/p$ and $\sigma=p-2$. Notice furthermore that 
we can write $r=\alpha p + q$ with $\alpha=(r-q)/p\ge 0$. The condition 
$r<(p-1)(q-1)/(p-2)$ ensures that $(p-2)\alpha \le \rho$. 
Let $0\le m\le p-2$ be arbitrary. We have $mr=m\alpha p+mq$ with
$0\le m\alpha\le (p-2)\alpha\le \rho$ and $0\le m\le \sigma=p-2$.
By Lemma \ref{phiflat} we then infer that $a_{pq}(mr)=1$. On invoking
Lemma \ref{verbinding} and Theorem \ref{upper} the proof of this case is then
completed.\\
{\tt Second case}. Assume that $q\equiv r\equiv 1({\rm mod~}p)$.\\
We claim that $r(p-2)\le (p-1)(q-1)-2$. By assumption we have
$r(p-2)< (p-1)(q-1)$. Suppose that $r(p-2)=(p-1)(q-1)-1$. By considering
this equation modulo $p$ we see that it is impossible and thus $r(p-2)\le (p-1)(q-1)-2$.
Note that $\sigma=0$ and $\rho=(p-1)(q-1)/p$. We can write $r=\alpha p + q$ with 
$\alpha=(r-q)/p\ge 0$. The condition 
$r(p-2)\le (p-1)(q-1)-2$ ensures that $(p-2)\alpha \le q-2-\rho$. 
Let $0\le m\le p-2$ be arbitrary. We have $1+mr=1+m\alpha p+mq$ with
$0\le m\alpha\le (p-2)\alpha\le q-2-\rho$ and $0\le m\le p-2-\sigma=p-2$.
By Lemma \ref{phiflat} we then infer that $a_{pq}(1+mr)=1$. On invoking
Lemma \ref{verbinding} and Theorem \ref{upper} the proof of this case is then
also completed.\qed\\

\noindent {\tt Remark}. (Y. Gallot.) The above result suggests perhaps that 
in case $n$ is of order at least two, $V_n$ is always of the form
$\{-a,-(a-1),\cdots,-1,0,1,\cdots,(a-1),a\}$ for some positive integer $a$.
However, this is not the case. The smallest $n$ for which $V_n$ is not
of this form is $n=23205 = 3\cdot 5\cdot 7\cdot 13\cdot 17$. 
Here the height is 13, but 12 (and -12) are not included in $V_n$.
Further examples (in order of appearance) are 46410 (height 13, $\pm 12$ not there), 
49335 (height  34, $\pm 33$ not found), 50505 (height 15,  $\pm 14$ not found). There
are also examples where a whole range values smaller than the height is not in $V_n$.

\subsection{The case where $p=3$}
In the case where $p=3$ we can always explicitly compute $V_{3qr}$ on
invoking Theorem \ref{ext}, Theorem \ref{extreme} and Lemma \ref{opstart}.
We obtain the following result.
\begin{Thm}
\label{drie}
Let $3<q<r$ be odd primes.\\ 
If $q\equiv 1({\rm mod~}3)$, $r\equiv 1({\rm mod~}3)$ 
and $r\le 2q-7$, then $V_{3qr}=\{-2,-1,0,1,2\}$. In
particular, $c_{3qr}(r+1)=2$ and $c_{3qr}(r+1+qr)=-2$.\\
If $q\equiv 2({\rm mod~}3)$, $r\equiv 2({\rm mod~}3)$ 
and $r\le 2q-3$, then $V_{3qr}=\{-2,-1,0,1,2\}$. In
particular, $c_{3qr}(r)=-2$ and $c_{3qr}(r+qr)=2$.\\
In the remaining cases $V_{3qr}=\{-1,0,1\}$ and then
$\Psi_{3qr}(x)$ is flat.
\end{Thm}
{\tt Remark}. The quoted results only give $r\le 2q-3$. Note, however, that if
$q\equiv r\equiv 1({\rm mod~}3)$ and $r\le 2q-3$, then $r\le 2q-7$.\\

\noindent We now infer some consequences of Theorem \ref{drie}. {}For this
we need the following generalisation of Bertrand's Postulate.
\begin{Lem}
\label{post}
If $q$ is any prime, then the interval $(q,2q-7]$ contains primes
$p_1$ and $p_2$ with $p_i\equiv i({\rm mod~}3)$.
\end{Lem}
{\it Proof}. Molsen \cite{mol}, cf.
Moree \cite{bertje}, has shown that for $x\ge 199$ the interval
$(x,{8\over 7}x]$ contains primes
$p_1$ and $p_2$ with $p_i\equiv i({\rm mod~}3)$. {}From this
the result follows after some easy computations. \qed

\begin{Thm} $~$\label{vier}\\
{\rm 1)} Let $r$ be any prime, then $\Psi_{15r}(x)$ and $\Psi_{21r}(x)$ are flat.\\
{\rm 2)} Let $q\ge 11$ be a prime. Then 
$\Psi_{3qr}$ is flat for all primes $r\ge 2q-1$. However, there is at least one prime $r$ such that
$\Psi_{3qr}(x)$ is non-flat.\\
{\rm 3)} Let $3<q<r$ be primes. {}For $k\le 16$ we have $|c_{3qr}(k)|\le 1$.
\end{Thm}
{\it Proof}. 1) An immediate consequence of Theorem \ref{drie} and  
Lemma \ref{kleiner3}.\\
2) A consequence of Theorem \ref{drie} and Lemma \ref{post}.\\
3) By part 1 and Theorem \ref{drie} we infer that the smallest $r$ for
which $V_{3qr}\ne \{-1,0,1\}$ is $r=17$. By Lemma \ref{flauw} the
proof is then completed. \qed

\subsection{Reciprocal polynomials of intermediary height}
A variation of the methods used to establish Theorem \ref{main} yields the following upper
bound for $h(\Psi_{pqr}(x))$. Sometimes this bound is actually optimal, for example for
the Chernick Carmichael numbers (see Lemma \ref{carm}).
\begin{Thm}
\label{sigma}
Let $\rho$ and $\sigma$ be the unique non-negative integers such that one
has
$(p-1)(q-1)=\rho p+\sigma q$. 
Put $\tau=(p-1)(q+r-1)$. If $qr>\tau$, then
the height of $\Psi_{pqr}(x)$ is at
most $\max\{\min(\rho+1,\sigma+1),\min(q-1-\rho,p-1-\sigma)\}$.
\end{Thm}
\begin{cor}
If either $q\equiv -2({\rm mod~}p)$ or $q\equiv 2({\rm mod~}p)$ and $q>p+2$,
then the height of 
$\Psi_{pqr}(x)$ is at most $(p+1)/2$.
\end{cor}
{\it Proof}. One easily checks that $qr>\tau$. We compute that
$$\sigma=\cases{{p-3\over 2} & if $q\equiv -2({\rm mod~}p)$;\cr
{p-1\over 2} & if $q\equiv 2({\rm mod~}p)$.}$$
{\it Proof of Theorem} \ref{sigma}. We have to show that $|c_{pqr}(k)|$ does
not exceed the bound stated. The
conditions of Lemma \ref{verbinding} are satisfied and by property (\ref{prop})
we may take $k\le \tau/2<(p-1)r$. 
Now choose $0\le m\le p-2$ such
that $mr\le k<(m+1)r$. By Lemma \ref{verbinding} we have
$$c_{pqr}(k)=-\sum_{v=0}^m a_{pq}(k-vr).$$
Let us consider the worst case where $m=p-2$ and a priori $|c_{pqr}(k)|\le p-1$.
We determine the maximum number of $v$ with 
$0\le v\le p-2$ for which $a_{pq}(k-vr)=1$. Let us
suppose that for $v_1,\ldots,v_t$ we have $a_{pq}(k-v_jr)=1$ and
hence, by Lemma \ref{phiflat}, we have
$$\cases{k-v_1r=i_1p+j_1q;\cr  k-v_2r=i_2p+j_2q; \cr \ldots  \cr k-v_tr=i_tp+j_tq,}$$
where each $j_m$ satisfies $0\le j_m\le \sigma$. Now if $t>\sigma+1$
two of the $j_m$ must be equal. On subtracting the corresponding equations
it would follow that $p|r$, a contradiction that shows that $t\le \sigma +1$. 
On using that $q\nmid r$, we likewise infer that $t\le \rho+1$.
We infer that
$c_{pqr}(k)\ge -\min(\rho+1,\sigma+1)$. Note that the same inequality actually holds for
all $k<(p-1)r$.\\
\indent We determine the maximum number of $w$ with 
$0\le w\le p-2$ for which $a_{pq}(k-wr)=-1$. Let us
suppose that for $w_1,\ldots,w_t$ we have $a_{pq}(k-w_jr)=1$ and
hence, by Lemma \ref{phiflat}, we have
$$\cases{k-w_1r=1+i_1p+j_1q;\cr  k-w_2r=1+i_2p+j_2q; \cr \ldots \cr k-w_tr=1+i_tp+j_tq,}$$
where each $j_m$ satisfies $0\le j_m\le p-2-\sigma$. Now if $t>p-1-\sigma$
two of the $j_m$ must be equal. On subtracting the corresponding equations
it would follow that $p|r$, a contradiction that shows that $t\le p-1-\sigma$. 
Likewise we infer that $t\le q-1-\rho$.
We infer that
$c_{pqr}(k)\le \min(q-1-\rho,p-1-\sigma)$. On combining this with 
$c_{pqr}(k)\ge  -\min(\rho+1,\sigma+1)$ we are done. \qed

\section{{}Further flatness results}
In this section we present some further (near) flatness results. 
\begin{Lem}
If $r>(p-1)(q-1)$, then $\Psi_{pqr}(x)$ is flat.
\end{Lem}
{\it Proof}. Note that if $f$ and $g$ are flat polynomials and $m>{\rm deg}(f)$, then
$f(x)g(x^m)$ is flat. By (\ref{pqr2}) we have $\Psi_{pqr}(x)=\Phi_{pq}(x)\Psi_{pq}(x^r)$.
The assumption on $r$ implies that $r>{\rm deg}(\Phi_{pq}(x))=(p-1)(q-1)$. Since
both $\Phi_{pq}(x)$ and $\Psi_{pq}(x)$ are flat, the result now follows. \qed\\

\noindent A variation of the latter proof making use of the identity
$\Psi_{pn}(x)=\Psi_n(x^p)\Phi_n(x)$ if $p\nmid n$ (this is part 3 of Lemma \ref{blup}),
yields the following lemma.
\begin{Lem} Let $p$ be a prime.
Let $h_1,h_2$ be the height
of  $\Phi_n(x)$, respectively $\Psi_n(x)$. If 
$p>\varphi(n)$, then $\Psi_{np}(x)$
is of height $h_1h_2$.
\end{Lem} 
Using this result we easily infer the following one.
\begin{Lem} Let $3<q<r<s$ be primes such that $s>2(q-1)(r-1)$. Then\\
{\rm 1)} $\Psi_{3qrs}(x)$ is of height at most $4$.\\
{\rm 2)} If $r\equiv q({\rm mod~}3)$ and $r\equiv \pm 1({\rm mod~}3q)$, then
$\Psi_{3qrs}(x)$ is flat.
\end{Lem}
{\it Proof}. 1) Beiter \cite{Beiter} has shown that $\Phi_{3qr}(x)$ is of height
at most 2. By Theorem \ref{drie} we know that also $\Psi_{3qr}(x)$
is of height at most 2. Now apply the previous lemma with
$n=3qr$ and $p=s$.\\
2) {}Follows from the previous lemma, Theorem \ref{drie} and the result due 
to Kaplan \cite[Theorem 1]{Kaplan} (who extended on earlier work by 
Bachman \cite{B3}) that
$\Phi_{3qr}(x)$ is flat if $r\equiv \pm 1({\rm mod~}3q)$.\qed\\

\noindent {\tt Remark}. Since $h(\Psi_{3\cdot 11\cdot 17\cdot 331}(x))=4$, we
see that the 4 above cannot be replaced by a smaller number.

Recall that smallest $n$ for which $\Phi_n(x)$ is non-flat
is $n=105$.
\begin{Lem} 
\label{105}
The smallest $n$ for which $\Psi_n(x)$ is non-flat
is $n=561$.
\end{Lem}
{\it Proof}. By computation one finds that $c_{561}(17)=-2$. 
By Lemma \ref{kleiner3} it suffices to check
that $\Psi_n(x)$ is flat for every odd 
squarefree $n\le 560$ with $\omega_1(n)\ge 3$.
This leaves us with the sets $${\cal
A}=\{105,165,195,231,255,273,285,345,357,399,435,465,483,555\},$$ and 
${\cal B}=\{385,429,455\}$, where the set ${\cal A}$ has
all its elements divisible by 15 or 21. On applying part 1 of Theorem
\ref{vier} we infer that $\Psi_n(x)$ is flat for every $n\in {\cal A}$.
By direct computation we find that $\Psi_{385}(x),\Psi_{429}(x)$ and 
$\Psi_{455}(x)$ are flat. \qed\\

\noindent Since 561 is the smallest Carmichael number and the 
smallest number $m$ for which $h(\Psi_m(x))>1$, one might wonder
whether perhaps $h(\Psi_C(x))>1$ for every Carmichael number $C$. The
answer is no, as the example $c=2821$ shows. However, for the
Chernick Carmichael numbers the answer turns out to be yes. In 1939 
Chernick proved that if $k\ge 0$ is such that $6k+1$, $12k+1$ and
$18k+1$ are all primes, then $C=(6k+1)(12k+1)(18k+1)$ is a Carmichael number.
Examples occur for $k=1,6,35,45,51,56,\ldots$.
\begin{Lem}
\label{carm}
If $C=(6k+1)(12k+1)(18k+1)$ is a Chernick Carmichael number, then
$c_{C}(24k+2)=-2$ and $h(\Psi_C(x))=2$.
\end{Lem}
{\it Proof}. Put $p=6k+1$, $q=12k+1$ and $r=18k+1$. We find $\rho=1$ and
$\sigma=p-2$. By Theorem \ref{sigma} we infer that $h(\Psi_C(x))\le 2$.
By Lemma \ref{phiflat} we have $a_C(2q)=1$ and $a_C(p)=1$.
Now $c_C(2q)=-a_C(2q)-a_C(2q-r)=-a_C(2q)-a_C(p)=-2$. 
Thus $c_C(2q)=c_C(24k+2)=-2$ and $h(\Psi_C(x))=2$.\qed

\section{Sizable coefficients}
The history of sizable coefficients goes back to Schur who in a letter in 1931 to
Landau (see e.g. E. Lehmer \cite{Emma}) proved that the $a_n(k)$ are unbounded. 
It is not difficult, see Suzuki \cite{Suzuki},
to adapt his argument so as to show that {\it every} integer shows up as a coefficient,
that is $\{a_n(k):n\ge 1,~k\ge 0\}=\mathbb Z$. Bungers \cite{Bungers}, in his Ph.D. thesis
proved that under the assumption that there are infinitely many twin primes, the
$a_n(k)$ are also unbounded if $n$ has at most three prime factors. E. Lehmer \cite{Emma} eliminated
the unproved assumption of the existence of infinitely twin primes from this. The strongest
result in this direction to date is due to Bachman, who proved a result 
(\cite[Theorem 1]{B2}), which implies that
$$\{a_{pqr}(k):3\le p<q<r~{\rm primes}\}=\mathbb Z.$$
\indent A minor variation of Suzuki's argument gives $\{c_n(k):n\ge 1,~k\ge 0\}=\mathbb Z$.
Since the next result is stronger, the details are left to the interested reader.
\begin{Thm}
We have $\{c_{pqr}(k):3\le p<q<r~{\rm primes}\}=\mathbb Z$.
\end{Thm}
{\it Proof}. By Dirichlet's theorem on arithmetic progressions for every prime $p$ 
there is a $q_0(p)$ such that for every $q>q_0(p)$ with $q\equiv \pm 1({\rm mod~}p)$,
there exists $r\equiv q({\rm mod~}p)$ with $q<r<(p-1)(q-1)/(p-2)$. The proof is
then completed on invoking
Theorem \ref{extreme}.\qed\\

\noindent In the table below (part of a much large table computed by Yves Gallot) the minimal $n$, $n_0$, such that $c_{n_0}(k)=m$ for
some $k$ is given. The third column gives the degree of $\Psi_{n_0}(x)$. The
fourth column gives the smallest $k$, $k_0$, for which $|c_{n_0}(k_0)|=m$.\\
\vfil\eject

\centerline{{\bf Table 1:} {\bf Minimal $n$ and $k$ with $|c_n(k)|=m$}}
\begin{center}
\begin{tabular}{|c|c|c|c|c|c|}
\hline
$m$    & $n_0$ & {\rm deg}$(\Psi_{n_0})$ & $k_0$  & $c_{n_0}(k_0)$\\
\hline
$1$    & $1$ & $0$ & $0$ & $+1$\\
\hline
$2$    & $561=3\cdot11\cdot 17$ & $241$ & $17$ & $-2$\\
\hline
$3$ & $1155=3\cdot 5\cdot 7\cdot 11$ & $675$ & $33$ & $-3$\\
\hline
$4$ & $2145=3\cdot 5\cdot 11\cdot 13$ & $1185$ & $44$ & $+4$\\
\hline
$5$  & $3795=3\cdot 5\cdot 11\cdot 23$ & $2035$ & $132$ & $-5$\\
\hline
$6$  & $5005=5\cdot 7\cdot 11\cdot 13$ & $2125$ & $201$ & $-6$\\
\hline
$7$  & $5005=5\cdot 7\cdot 11\cdot 13$ & $2125$ & $310$ & $-7$\\
\hline
$8$  & $8645=5\cdot 7\cdot 13\cdot 19$ & $3461$ & $227$ & $-8$\\
\hline
$9$  & $8645=5\cdot 7\cdot 13\cdot 19$ & $3461$ & $240$ & $+9$\\
\hline
$10$  & $11305=5\cdot 7\cdot 17\cdot 19$ & $4393$ & $240$ & $-10$\\
\hline
$11$  & $11305=5\cdot 7\cdot 17\cdot 19$ & $4393$ & $306$ & $+11$\\
\hline
\end{tabular}
\end{center}
{}For $m=10,\ldots,21$ it turns out that $n_0=11305$.\\

\noindent {\bf Acknowledgement}. As in many papers before this one, impressive computational 
assistance was provided by Yves Gallot. I am deeply grateful for all the work
he did for me. Paul Tegelaar proofread an earlier version and Gennady
Bachman updated me on recent results on the flatness of $\Phi_{3qr}(x)$.\\
\indent In the context of assigning student projects to 
Alexander Bridi, Andreas Decker, Patrizia Dressler, Silke Glas and Thorge Jensen (and
earlier Christine Jost and Janina M\"uttel),
I considered cyclotomic polynomials. I thank them all for their enthusiasm and
cheerful presence which added to the already pleasant atmosphere at the
Max-Planck-Institut f\"ur Mathematik.\\
\indent My belief that $561$ is a boring number was expelled by Don Zagier who remarked
that, on the contrary, it is interesting since it is the smallest Carmichael  number. I
thank him for this and several other helpful remarks.

\medskip\noindent {\footnotesize Max-Planck-Institut f\"ur Mathematik,\\
Vivatsgasse 7, D-53111 Bonn, Germany.\\
e-mail: {\tt moree@mpim-bonn.mpg.de}}
\end{document}